\documentclass{amsart}

\usepackage[T1]{fontenc}
\usepackage{lmodern}
\usepackage{amsmath,amssymb,amsthm,mathtools}
\usepackage{booktabs,tabularx,array}
\usepackage{microtype}
\usepackage[hidelinks]{hyperref}

\newif\ifanon
\anonfalse

\ifanon
\hypersetup{pdftitle={A Complete Classification of Maximal Hyperelliptic Curves Related to Chebyshev Polynomials},
pdfauthor={}}
\else
\hypersetup{pdftitle={A Complete Classification of Maximal Hyperelliptic Curves Related to Chebyshev Polynomials},
pdfauthor={Saeed Tafazolian and Jaap Top}}
\fi
\hypersetup{pdfsubject={Maximal curves, Chebyshev polynomials, and complex multiplication},
pdfkeywords={Chebyshev polynomials, maximal curves, hyperelliptic curves, complex multiplication, Dirichlet L-functions}}

\newtheorem{theorem}{Theorem}[section]
\newtheorem{proposition}[theorem]{Proposition}
\newtheorem{lemma}[theorem]{Lemma}
\newtheorem{corollary}[theorem]{Corollary}
\theoremstyle{remark}
\newtheorem{remark}[theorem]{Remark}

\newcommand{\F}{\mathbb F}
\newcommand{\Q}{\mathbb Q}
\newcommand{\C}{\mathbb C}
\newcommand{\Z}{\mathbb Z}
\newcommand{\Hh}{\mathcal H}
\newcommand{\End}{\operatorname{End}}
\newcommand{\Gal}{\operatorname{Gal}}
\newcommand{\Jac}{\operatorname{Jac}}
\newcommand{\ord}{\operatorname{ord}}

\title[Maximal Chebyshev hyperelliptic curves]{A Complete Classification of Maximal Hyperelliptic Curves Related to Chebyshev Polynomials}
\ifanon\else
\author{Saeed Tafazolian}
\address{Departamento de Matem\'atica, IMECC, Universidade Estadual de Campinas,
Campinas, SP, Brazil}
\email{saeed@unicamp.br}
\author{Jaap Top}
\address{Bernoulli Institute for Mathematics, Computer Science, and Artificial
Intelligence, University of Groningen, Nijenborgh 9, 9747 AG Groningen,
The Netherlands}
\email{j.top@rug.nl}
\fi
\date{}

\begin{document}

\begin{abstract}
Let $\varphi_d$ be the Chebyshev polynomial with
$\varphi_d(t+t^{-1})=t^d+t^{-d}$, and let $\Hh_d$ be the hyperelliptic
curve $y^2=\varphi_d(x)$. We determine all pairs $(d,p)$, with $p\nmid2d$
prime, for which $\Hh_d$ is maximal over some finite field of
characteristic $p$, together with all such fields. For
$p\equiv1\pmod4$ only partial results were known. In this case, for odd $d$ we show that the Jacobian of $\Hh_d$ is
supersingular if and only if $\Hh_d$ is maximal over a finite extension of
$\F_p$, and that this happens precisely when $d=\ell^n$ with
$\ell\equiv1\pmod4$ prime and $p$ a generator of $(\Z/\ell^n\Z)^\times$;
we also compute the Frobenius polynomial. Necessity comes from the
Shimura--Taniyama formula together with the non-vanishing of
the value at $0$ of certain Dirichlet L-series and information about the action of
Frobenius on $2$-torsion. The case of even $d$ and $p\equiv1\pmod 4$ follows from the odd case, by
a short reduction.
\end{abstract}

\maketitle

\section{Introduction}\label{sec:intro}

A curve of genus $g$ over $\F_q$ is maximal if it attains the Hasse--Weil
bound $q+1+2g\sqrt q$. Maximal curves of positive genus are supersingular,
they exist only over fields of square order, and they have been studied
extensively, in
part because of their role in coding theory \cite{Stichtenoth}. This paper
completes the answer to the following question: for which $d$ and $p$ is
the hyperelliptic curve defined by a Chebyshev polynomial maximal over some
finite field of characteristic $p$?

For an integer $d\geq2$, let $\varphi_d\in\Z[x]$ be the Chebyshev
polynomial characterized by
\begin{equation}\label{eq:chebdef}
 \varphi_d(t+t^{-1})=t^d+t^{-d}.
\end{equation}
Thus $\varphi_d(x)=D_d(x,1)=2T_d(x/2)$ in the usual Dickson/Chebyshev
notation. If $p\nmid2d$ is prime, we write $\Hh_d/\F_p$ for the smooth
projective curve with affine equation $y^2=\varphi_d(x)$. Its genus is
$\lfloor(d-1)/2\rfloor$.

The question was considered in \cite{TT19,DTT25,TT26}. For
$p\equiv3\pmod4$ the answer is known: for odd $d$ it is
\cite[Proposition~4.1.1]{TT26}, and for even $d$ it follows from
\cite{TT19,DTT25}. The case $p\equiv1\pmod4$ behaves differently. For odd
$d$ one has $\varphi_d(x)=xg(x^2)$ with $g$ separable and $g(0)\neq0$.
For curves of this shape, \cite[Proposition~2.3]{TT19} asserted that a
supersingular Jacobian forces $p\equiv3\pmod4$; this assertion was subsequently corrected in \cite{TTcorr}, and Theorem~\ref{thm:classification} describes
exactly when it fails for the curves $\Hh_d$. In \cite[Theorem~3.4.1]{TT26}
it is shown that $\Hh_\ell$ is maximal over a finite extension of $\F_p$
when $\ell\equiv1\pmod4$ is prime and $p\equiv1\pmod4$ is a primitive
root modulo $\ell$. The same paper constructs CM factors for degrees
$\ell_1\ell_2$ and $\ell^n$, and its Questions~3.4.8, 4.2.4 and~4.2.9 ask
for the converse and for the corresponding statement in general odd
degree. We answer these questions.

\begin{theorem}\label{thm:classification}
Let $d\geq3$ be odd and let $p\nmid2d$ be a prime with
$p\equiv1\pmod4$. The following conditions are equivalent:
\begin{enumerate}
\item $\Jac(\Hh_d)$ is supersingular;
\item $\Hh_d$ is maximal over some finite extension of $\F_p$;
\item $d=\ell^n$ for a prime $\ell\equiv1\pmod4$ and an integer
$n\geq1$, and $p$ generates $(\Z/\ell^n\Z)^\times$.
\end{enumerate}
\end{theorem}

We write $P_{C,p}(X)$ for the characteristic polynomial of the $p$-power
Frobenius of a curve $C/\F_p$, or of its Jacobian.

\begin{theorem}\label{thm:frobenius}
Let $\ell\equiv1\pmod4$ and $p\equiv1\pmod4$ be distinct primes, let
$n\geq1$, and suppose that $p$ generates $(\Z/\ell^n\Z)^\times$. Put
\[
 m_j=\ell^{j-1}(\ell-1)\qquad(1\leq j\leq n).
\]
Then
\begin{equation}\label{eq:mainpoly}
 P_{\Hh_{\ell^n},p}(X)
 =\prod_{j=1}^n\left(X^{m_j}+p^{m_j/2}\right).
\end{equation}
For a positive integer $r$, the curve $\Hh_{\ell^n}$ is maximal over
$\F_{p^r}$ if and only if $r=m_nb$ with $b$ positive and odd. It is
minimal over $\F_{p^r}$ if and only if $r=2m_nb$ with $b\geq1$.
\end{theorem}

For fixed $\ell\equiv1\pmod4$ and $n\geq1$, there are infinitely many
primes satisfying the hypotheses of Theorem~\ref{thm:frobenius}: prescribe
a primitive root modulo $\ell^n$ and the residue $1$ modulo $4$, and apply
the Chinese remainder theorem and Dirichlet's theorem.

The odd-degree result also gives the remaining even-degree cases; the
proof is in Section~\ref{sec:even}.

\begin{theorem}\label{thm:even-p1}
Let $p\equiv1\pmod4$ be prime, let $d\geq4$ be even, and assume
$p\nmid2d$. Then $\Hh_d$ is maximal over some finite extension of
$\F_p$ if and only if one of the following holds:
\begin{enumerate}
\item $d=4$ and $p\equiv5\pmod8$;
\item $d=2\ell^a$, where $\ell\equiv1\pmod4$ is prime, $a\geq1$, and
$p$ generates $(\Z/\ell^a\Z)^\times$.
\end{enumerate}
In case \emph{(1)}, the maximal fields are precisely $\F_{p^{2b}}$ with
$b$ positive and odd. In case \emph{(2)}, if
$M=\ell^{a-1}(\ell-1)$, they are precisely $\F_{p^{Mb}}$ with $b$
positive and odd.
\end{theorem}

Together with \cite[Proposition~4.1.1]{TT26} and the even-degree results
of \cite{TT19,DTT25}, this gives a classification for all $d\geq3$,
collected in Table~\ref{tab:complete} and Corollary~\ref{cor:all-d}.

\subsection*{Method}
For the CM factors of \cite{TT26}, the Shimura--Taniyama formula turns
supersingularity into a balance condition on cosets of a decomposition
group, and finite Fourier analysis turns this into the vanishing of
$L(0,\chi\chi_4)$ for the even characters $\chi$ trivial on that group.
Non-vanishing of these values gives the necessity of~(3) in
Theorem~\ref{thm:classification}. Supersingularity does not determine the
Frobenius polynomial; for Theorem~\ref{thm:frobenius} the remaining sign
is fixed by the action of Frobenius on geometric $2$-torsion.

\subsection*{Relation to the cover and to Gauss sums}
The curve $\Hh_d$ is the quotient of $X_d:w^2=t^{2d+1}+t$ by an
involution. For this cover, Kodama and Washio \cite{KodamaWashio} proved
that supersingularity, minimality over an extension and maximality over an
extension are equivalent; see \cite[Remark~4.1.2]{TT26}. When
$p\equiv1\pmod4$, however, $\Jac(X_d)$ is never supersingular
(Remark~\ref{rem:cover}), so Theorem~\ref{thm:classification} is not
inherited from the cover. Moreover, with $z=-t^{2d}$ the curve $X_d$ is
birational to $w^{4d}=-z(1-z)^{2d}$, a quotient of a Fermat curve, so its Frobenius eigenvalues over a large enough field are Jacobi sums
\cite{Weil}. In the situation of Theorem~\ref{thm:frobenius},
$\langle p\rangle\leq(\Z/4\ell^n\Z)^\times$ has index $2$ and does not
contain $-1$, so the semiprimitive criterion for Fermat quotients
\cite{ShiodaKatsura} does not apply and the Gauss sums involved are of
index~$2$ \cite{YangXia}. The argument below avoids their evaluation and
works directly over $\F_p$.

\section{Preliminaries}\label{sec:prelim}

We use $\phi(d)$ for Euler's totient function, to distinguish it from
the Chebyshev polynomial $\varphi_d$. Throughout the paper, a curve is
understood to be smooth, projective and geometrically connected.

\subsection{Frobenius and isogeny factors}
For an abelian variety $A/\F_p$ of dimension $g$ and a prime
$\lambda\neq p$, write $T_\lambda(A)$ for the associated Tate module and put
\[
 V_\lambda(A)=T_\lambda(A)\otimes_{\Z_\lambda}\Q_\lambda,
 \qquad P_{A,p}(X)=\det(X-F_p\mid V_\lambda(A)).
\]
Here $F_p$ is the $p$-power Frobenius endomorphism. The polynomial
$P_{A,p}\in\Z[X]$ is independent of $\lambda$, has degree $2g$ and
constant term $p^g$, and all its complex roots have absolute value
$p^{1/2}$. Frobenius acts semisimply on $V_\lambda(A)$, and the action of
endomorphisms of $A$ on $V_\lambda(A)$ is faithful
\cite{Tate,MilneAV}.

For a curve $C/\F_p$ of genus $g$, write
$P_{C,p}=P_{\Jac(C),p}=\prod_{i=1}^{2g}(X-\alpha_i)$. Then
\begin{equation}\label{eq:pointcount}
 \#C(\F_{p^r})=p^r+1-\sum_{i=1}^{2g}\alpha_i^r.
\end{equation}
We call $C$ maximal, respectively minimal, over $\F_{p^r}$ if it
attains the upper, respectively lower, Hasse--Weil bound. For $g>0$
this requires $r$ to be even, and is equivalent to
\begin{equation}\label{eq:extremal}
 \alpha_i^r=-p^{r/2},\quad\text{respectively}\quad
 \alpha_i^r=p^{r/2}\qquad(1\leq i\leq2g).
\end{equation}
We use the standard characterization of supersingularity by the condition
that all slopes of the Newton polygon at $p$
of the polynomial $P_{C,p}$ are equal to $1/2$; see, for example, \cite{Oort}. Formula
\eqref{eq:extremal} immediately implies that a maximal or minimal curve
has supersingular Jacobian. Supersingularity alone does not in general
imply maximality over an extension. For instance, a supersingular elliptic
curve over $\F_{p^2}$ whose $p^2$-Frobenius is $[p]$ is minimal over every
finite extension of $\F_{p^2}$ and maximal over none.

If $f:C\to C'$ is a finite morphism of degree $c$ defined over $\F_p$,
then $f_*f^*=[c]$, so $\Jac(C')$ is an isogeny factor of $\Jac(C)$ and
supersingularity and maximality pass to $C'$
\cite{Stichtenoth,MilneAV}.

\subsection{Chebyshev subcovers}
The defining identity \eqref{eq:chebdef} gives
\begin{equation}\label{eq:composition}
 \varphi_{ab}=\varphi_a\circ\varphi_b=\varphi_b\circ\varphi_a.
\end{equation}
For $e\mid d$ there is consequently a degree-$d/e$ map
\begin{equation}\label{eq:covering}
 \pi_{d,e}:\Hh_d\longrightarrow\Hh_e,
 \qquad(x,y)\longmapsto(\varphi_{d/e}(x),y).
\end{equation}

\subsection{CM slopes}
Let $A$ be an abelian variety of dimension $h$ over a number field,
with an embedding $K\hookrightarrow\End^0(A_{\overline\Q})$ for a
CM field $K$ of degree $2h$; for simplicity we
moreover assume that $A$ is in fact defined over $K$. We read the CM type $\Phi$ from the
pullback action of $K$ on $H^0(A,\Omega_A^1)$. This gives the same
eigencharacters as the differential action on $\operatorname{Lie}(A)$;
see \cite[Chapter~I, \S3]{MilneCM}.

Fix embeddings $\overline{\Q}\hookrightarrow\C$ and
$\overline{\Q}\hookrightarrow\overline{\Q}_p$. Suppose that $K/\Q$ is
Galois, and use the first embedding to identify $\Phi$ with a subset of
$\Gal(K/\Q)$; the second fixes the prime of $K$ above $p$ and hence its
decomposition group $D$. At a prime of good reduction, the CM slope formula
gives that the set of slopes of the corresponding Newton polygon of the characteristic polynomial
of Frobenius equals the set
\begin{equation}\label{eq:CMslope}
\left\{ \frac{|\Phi\cap D\tau|}{|D|}
 \quad : \quad\tau\in \Gal(K/\Q)\right\}.
\end{equation}
This is the Shimura--Taniyama formula; see \cite{ShimuraTaniyama,TateBourbaki}
and, for the above formulation, \cite[Theorem~3.1]{Blake}. The formula also applies when the CM action is defined only after a
finite extension, since Newton slopes do not change under finite
extensions of the residue field. In particular, the reduction is
supersingular precisely when every intersection 
$\Phi \cap D\tau$ in \eqref{eq:CMslope}
has cardinality $|D|/2$.

For odd $d>1$, put
\begin{equation}\label{eq:CMfields}
 K_d=\Q(i,\zeta_d+\zeta_d^{-1}),\qquad
 M_d=\Q(i,\zeta_d),\qquad U_d=(\Z/d\Z)^\times/\{\pm1\}.
\end{equation}
Then $[K_d:\Q]=\phi(d)$ and
$\Gal(K_d/\Q)=\{\pm1\}\times U_d$, where $(\delta,\pm a)$ in the latter group is identified with
the automorphism of $K_d$ that sends
$i$ to $\delta i$ and $\zeta_d+\zeta_d^{-1}$ to
$\zeta_d^a+\zeta_d^{-a}$. At $p\nmid2d$ with $p\equiv1\pmod4$,
the decomposition group is
\begin{equation}\label{eq:decomp}
 D=\{1\}\times H\quad\text{where}\quad H=\langle\pm p\rangle\leq U_d.
\end{equation}
For $u=\pm a\in U_d$, let $\langle a\rangle$ be its unique
representative in $[1,(d-1)/2]$, and define
\begin{equation}\label{eq:epsilon}
 \epsilon_d(u)=(-1)^{\langle a\rangle+1}.
\end{equation}
This is a function on $U_d$, not in general a character.

\subsection{Characters and special values}
Characters of $U_d$ are identified with even Dirichlet characters modulo~$d$, 
extended by zero on nonunits. Let $\chi_4$ be the nontrivial
character modulo $4$. For a character $\eta$ modulo $N$,
\begin{equation}\label{eq:L0}
 L(0,\eta)=-\frac1N\sum_{b=1}^N b\eta(b);
\end{equation}
this follows from the expression of $L(s,\eta)=N^{-s}\sum_{r=1}^N \eta(r)\zeta(s,\frac{r}{N})$ where 
$\zeta(s,x)$ denotes the Hurwitz zeta
function defined (for $\text{Re}(s)>1$) by
$\zeta(s,x)=\sum\sum_{n=0}^\infty (n+x)^{-s}$;
its meromorphic continuation
satisfies $\zeta(0,x)=\tfrac12-x$. If
$\eta$ is primitive and odd, the functional equation gives $L(0,\eta)$ as
a nonzero multiple of $L(1,\overline\eta)$, and Dirichlet's theorem yields
\begin{equation}\label{eq:Lnonzero}
 L(0,\eta)\neq0;
\end{equation}
see \cite[Chapter~4]{Washington}, \cite[Chapter~4]{Davenport} and
\cite[\S25.15]{DLMF}.
We shall also use finite Fourier inversion: a function on a finite abelian
group vanishes if and only if all its character sums vanish.

\subsection{\texorpdfstring{Hyperelliptic $2$-torsion}{Hyperelliptic 2-torsion}}
Let $C:y^2=f(x)$ over a field of characteristic different from $2$, with
$f$ squarefree of degree $2g+1$. For each root $\alpha$, the class
$D_\alpha=[(\alpha,0)-\infty]$ lies in $\Jac(C)[2]$. These classes
generate the geometric $2$-torsion, with the single relation
\begin{equation}\label{eq:torsionrelation}
 \sum_{f(\alpha)=0}D_\alpha=0;
\end{equation}
see \cite[Chapter~IIIa]{MumfordTheta}. Thus the $2$-torsion is the
permutation module on the finite roots, modulo the sum of its basis
vectors. This description respects the Galois action.

\section{Primitive CM factors}\label{sec:CM}

The factors below are those of \cite[Lemmas~4.2.1 and~4.2.5]{TT26};
the prime-degree construction goes back to \cite[Proposition~4]{TTV}.
We give the construction with its Galois equivariance, which will be
needed to determine Frobenius. All factors in this section are initially
defined over $\Q$; the same notation will denote their good reductions.

\subsection{Specialization}

\begin{lemma}\label{lem:specialization}
Let $A/\Q$ have good reduction $A_0/\F_p$ at $p$. Let $E\subseteq L$ be finite
Galois extensions of $\Q$, with $p$ unramified in $L$. Suppose that
$\iota:E\hookrightarrow\End^0(A_L)$ satisfies
\[
 {}^g\iota(u)=\iota(g|_E(u))\qquad(g\in\Gal(L/\Q),\ u\in E).
\]
Fix a prime $\mathfrak P$ of $L$ above $p$, and let $\sigma$ be the
restriction to $E$ of its arithmetic Frobenius. Specialization induces
an embedding $\iota_0:E\hookrightarrow\End^0(A_{0,\overline{\F}_p})$
satisfying
\begin{equation}\label{eq:semilinear}
 F_p\iota_0(u)=\iota_0(\sigma(u))F_p.
\end{equation}
Rational idempotents specialize to the corresponding isogeny factors.
\end{lemma}

\begin{proof}
Endomorphisms extend to the abelian scheme over
$\mathcal O_{L,\mathfrak P}$, and specialization is injective because it
is compatible with the faithful action on prime-to-$p$ Tate modules
\cite{SerreTate,MilneAV}. For the reduction $f$ of an endomorphism,
$F_pf=f^{(p)}F_p$, and since $p$ is unramified in $L$, the twist
$f^{(p)}$ is the reduction of the conjugate of $f$ by the arithmetic
Frobenius at $\mathfrak P$. The equivariance hypothesis gives
\eqref{eq:semilinear}. Specialization is a homomorphism of rational
endomorphism algebras, so it preserves idempotents.
\end{proof}

\subsection{Differentials and isogeny factors}
Here we briefly recall and extend some ideas from \cite[\S4.2]{TT26}. For odd $d$, consider the auxiliary curve
\[
 X_d:\quad z^2=t^d+t^{-d}.
\]
The substitution $w=t^{(d+1)/2}z$ gives the genus-$d$ model
$w^2=t(t^{2d}+1)$. Define
\[
 s(t,z)=(t^{-1},z),\qquad
 \rho(t,z)=(\zeta_dt,z),\qquad
 \jmath(t,z)=(-t,iz).
\]
The quotient by $s$ is
\[
 q_d:X_d\longrightarrow\Hh_d,\qquad(t,z)\longmapsto(t+t^{-1},z).
\]

\begin{lemma}\label{lem:differentials}
The forms
\[
 \omega_a=(t^{a-1}-t^{-a-1})\frac{dt}{z},
 \qquad 1\leq a\leq(d-1)/2,
\]
are a basis of $q_d^*H^0(\Hh_d,\Omega^1)$. Under $X_d\to X_e$,
$(t,z)\mapsto(t^c,z)$, where $d=ce$, one has
\begin{equation}\label{eq:pullback}
 \omega_b\longmapsto c\omega_{cb}.
\end{equation}
Consequently:
\begin{enumerate}
\item for $d=\ell^j$ and $e=\ell^{j-1}$, the old space is indexed by
$\ell\mid a$, and its complement is indexed by $(a,d)=1$;
\item for $d=\ell_1\ell_2$ with distinct odd primes, the two old spaces
are indexed by $\ell_2\mid a$ and $\ell_1\mid a$. They are disjoint,
and their complement is indexed by $(a,d)=1$.
\end{enumerate}
\end{lemma}

\begin{proof}
On the model in $(t,w)$,
\[
 \omega_a=
 \left(t^{(d+2a-1)/2}-t^{(d-2a-1)/2}\right)\frac{dt}{w}.
\]
Both exponents lie between $0$ and $d-1$, so the forms are regular.
They are $s$-invariant and linearly independent. Since $q_d$ is the
separable quotient by $\langle s\rangle$, pullback identifies
$H^0(\Hh_d,\Omega^1)$ with the $s$-invariant regular differentials on $X_d$
(see, for example, \cite{Stichtenoth}); their number is $g(\Hh_d)$. Hence
they form the asserted basis. Substitution of
$u=t^c$ gives \eqref{eq:pullback}. The two descriptions of the old
spaces follow. In the second case, no index in the prescribed interval
is divisible by $\ell_1\ell_2$, so the two old spaces are disjoint.
\end{proof}

In particular, at prime-power level the old indices are multiples of
$\ell$, rather than $\ell^{j-1}$. Formula~\eqref{eq:pullback} gives the
corresponding index calculation in the proof of
\cite[Lemma~4.2.5]{TT26}. The CM type and the statements of that paper
remain unchanged.

Write $J_d=\Jac(\Hh_d)$. For an odd prime $\ell$, set $B_1=J_\ell$.
For $j\geq2$, let
\[
 f_j=\pi_{\ell^j,\ell^{j-1}},\qquad B_j=(\ker f_{j,*})^0.
\]

\begin{lemma}\label{lem:split}
For $j\geq2$, the map
\[
 J_{\ell^{j-1}}\times B_j\longrightarrow J_{\ell^j},
 \qquad(a,b)\longmapsto f_j^*a+b
\]
is a $\Q$-isogeny of $\ell$-power degree. For every $j\geq1$,
\[
 \dim B_j=\frac{\ell^{j-1}(\ell-1)}2.
\]
At $p\nmid2\ell$, the factors have good reduction and
\begin{equation}\label{eq:torsionsplit}
 J_{\ell^j}[2]\simeq J_{\ell^{j-1}}[2]\oplus B_j[2]
 \qquad(j\geq2)
\end{equation}
as geometric Frobenius modules.
\end{lemma}

\begin{proof}
The identity $f_{j,*}f_j^*=[\ell]$ implies surjectivity of $f_{j,*}$
and gives the dimension. If $f_j^*a+b=0$, applying $f_{j,*}$ gives
$[\ell]a=0$, and $b$ is determined by $a$. Thus the kernel embeds in
$J_{\ell^{j-1}}[\ell]$. It is finite of $\ell$-power order, and the
source and target have equal dimension, proving the isogeny assertion.

The Jacobian $J_{\ell^j}$ has good reduction at $p\nmid2\ell$.
The rational Tate module of its $\Q$-isogeny factor $B_j$ is a
Galois-stable direct summand, hence unramified there. The
N\'eron--Ogg--Shafarevich criterion gives good reduction of $B_j$
\cite{SerreTate}. Specializing the isogeny, whose degree is odd and
prime to $p$, gives \eqref{eq:torsionsplit}.
\end{proof}

For $d=\ell_1\ell_2$ with distinct odd primes, define
\begin{equation}\label{eq:Ad}
 A_d=\left(\ker\bigl(
 \pi_{d,\ell_1,*},\pi_{d,\ell_2,*}\bigr)\right)^0.
\end{equation}
Lemma~\ref{lem:differentials} shows that the combined norm map is
surjective: its pullback on differentials has image the direct sum of
the two old spaces. Thus
\begin{equation}\label{eq:twoprimesplit}
 J_d\sim_{\Q}J_{\ell_1}\times J_{\ell_2}\times A_d,
 \qquad\dim A_d=\frac{\phi(d)}2.
\end{equation}
Here the isogeny follows from Poincar\'e complete reducibility over $\Q$;
see \cite{MilneAV}.
The same good-reduction argument applies to $A_d$ at $p\nmid2d$.

\subsection{The CM action}

\begin{proposition}\label{prop:CM}
Let $d=\ell^j$ or $d=\ell_1\ell_2$ with distinct odd primes, and let
$A=B_j$ or $A=A_d$, respectively. There is an embedding
\[
 \iota_d:K_d\hookrightarrow\End^0(A_{M_d}),\qquad
 [K_d:\Q]=2\dim A,
\]
with CM type
\begin{equation}\label{eq:CMtype}
 \Phi_d=\{((-1)^{a+1},\pm a):
 1\leq a\leq(d-1)/2,\ (a,d)=1\}
       =\{(\epsilon_d(u),u):u\in U_d\}.
\end{equation}
The embedding is Galois-equivariant. At every $p\nmid2d$, its
specialization is faithful. If $p\equiv1\pmod4$, then
\begin{equation}\label{eq:CMtwist}
 F_p\iota_d(u)=\iota_d(\sigma_p(u))F_p,
 \qquad\sigma_p=(1,\pm p)\in\Gal(K_d/\Q),
\end{equation}
where $\iota_d$ also denotes the specialized action.
\end{proposition}

\begin{proof}
Since $s\rho s=\rho^{-1}$, the operator $\rho+\rho^{-1}$ preserves
the $s$-invariant part of $\Jac(X_d)$. On $J_d$ define the rational
correspondence
\[
 T_d=\tfrac12q_{d,*}(\rho+\rho^{-1})q_d^*.
\]
It is defined over $K_d^+=\Q(\zeta_d+\zeta_d^{-1})$ (and not in general
over $\Q$).
The automorphism $\jmath$ commutes with $s$ and descends to
$I_d:(x,y)\mapsto(-x,iy)$ on $\Hh_d$. The induced endomorphisms
$T_d$ and $I_d$ commute, and $I_d^2=[-1]$. Using
$q_{d,*}q_d^*=[2]$ and the basis above, we obtain
\begin{equation}\label{eq:CMeigen}
 T_d^*\omega_a=(\zeta_d^a+\zeta_d^{-a})\omega_a,
 \qquad I_d^*\omega_a=(-1)^{a+1}i\omega_a.
\end{equation}

Let $F$ denote the norm map defining $A$, with $F=0$ when $j=1$,
and let $\nu:A\hookrightarrow J_d$ be the inclusion. The old
differential space $\operatorname{im}F^*$ is preserved by both operators.
For $R=T_d$ or $I_d$, the pullback of $FR\nu$ is
$\nu^*R^*F^*=0$, because $\nu^*$ annihilates the old space.
In characteristic zero a homomorphism of abelian varieties with zero
differential is zero; for rational endomorphisms one first clears
denominators. Thus both operators restrict to $A$. The primitive
span, indexed by $(a,d)=1$, maps isomorphically to $H^0(A,\Omega_A^1)$.

On this space, the eigenvalues of $T_d$ are the distinct conjugates
of $\zeta_d+\zeta_d^{-1}$. The faithful action on differentials
therefore identifies $\Q[T_d|_A]$ with the real cyclotomic field
$K_d^+=\Q(\zeta_d+\zeta_d^{-1})$. Since $X^2+1$ is irreducible over
this totally real field, adjoining the commuting $I_d$ gives
\[
 \Q[T_d|_A,I_d|_A]\simeq K_d^+[X]/(X^2+1)=K_d.
\]
Its degree is $\phi(d)=2\dim A$, and \eqref{eq:CMeigen} gives
\eqref{eq:CMtype}.

If $g(\zeta_d)=\zeta_d^r$ and $g(i)=\delta i$, then
\[
 {}^gT_d=\tfrac12q_{d,*}(\rho^r+\rho^{-r})q_d^*,
 \qquad {}^gI_d=\delta I_d.
\]
The first operator acts on $\omega_a$ by
$\zeta_d^{ra}+\zeta_d^{-ra}$, so it represents the natural image
of the real cyclotomic generator under $g$. This proves equivariance
on both generators. The extension $M_d/\Q$ is unramified at
$p\nmid2d$, and the factors have good reduction there. Applying
Lemma~\ref{lem:specialization} gives faithful specialization and
\eqref{eq:CMtwist}.
\end{proof}

\section{The Fourier obstruction and necessity}\label{sec:Fourier}

\begin{proposition}\label{prop:balance}
Let $p\nmid2d$ be a prime with $p\equiv1\pmod4$, and let $A$ be an
abelian variety over $\Q$ with good reduction at $p$ and with CM by $K_d$
of type \eqref{eq:CMtype}, as in Proposition~\ref{prop:CM}. Put
$H=\langle\pm p\rangle\leq U_d$. The reduction of $A$ is supersingular if
and only if
\begin{equation}\label{eq:balance}
 \sum_{u\in aH}\epsilon_d(u)=0\qquad(a\in U_d).
\end{equation}
Equivalently, for every even Dirichlet character $\chi$ modulo $d$
that is trivial on $H$,
\begin{equation}\label{eq:charsum}
 S_d(\chi):=\sum_{a=1}^{(d-1)/2}(-1)^{a+1}\chi(a)=0.
\end{equation}
\end{proposition}

\begin{proof}
By \eqref{eq:decomp}, the cosets with first coordinate $+1$ and $-1$
count the odd and even distinguished representatives, respectively.
Their slopes are
\[
 \frac12\pm\frac1{2|H|}\sum_{u\in aH}\epsilon_d(u),
\]
proving the first assertion. Define $F(aH)=\sum_{u\in aH}\epsilon_d(u)$
on $U_d/H$. For a character $\chi$ of this quotient,
\[
 \sum_{aH\in U_d/H}F(aH)\overline{\chi(a)}
 =\sum_{u\in U_d}\epsilon_d(u)\overline{\chi(u)}.
\]
Fourier inversion, followed by replacing $\chi$ with $\overline\chi$,
gives \eqref{eq:charsum}.
\end{proof}

\begin{lemma}\label{lem:Lidentity}
Let $d>1$ be odd and $\chi$ an even Dirichlet character modulo $d$.
Then
\begin{equation}\label{eq:Lidentity}
 S_d(\chi)=-\frac{\chi_4(d)}{\chi(2)}L(0,\chi\chi_4),
\end{equation}
where the product is regarded as a character modulo $4d$.
If $\chi$ is primitive and nontrivial, then $\eta=\chi\chi_4$ is
primitive and odd, $L(0,\eta)\neq0$, and hence $S_d(\chi)\neq0$.
\end{lemma}

\begin{proof}
Put $\eta=\chi\chi_4$. It is odd and nontrivial. Since
$\eta(b+2d)=-\eta(b)$, pairing terms in \eqref{eq:L0} gives
\[
 L(0,\eta)
 =-\frac1{4d}\sum_{b=1}^{4d}b\eta(b)
 =\frac12\sum_{b=1}^{2d}\eta(b).
\]
Moreover $\eta(2d-b)=\eta(b)$, and the terms at $d$ and $2d$ vanish.
Hence $L(0,\eta)=\sum_{b=1}^{d-1}\eta(b)$. Only odd $b$ contribute;
write $b=d-2a$, with $1\leq a\leq(d-1)/2$. Evenness of $\chi$ and
the identity $\chi_4(d-2a)=\chi_4(d)(-1)^a$ give
\[
 L(0,\eta)
 =\chi(2)\chi_4(d)\sum_{a=1}^{(d-1)/2}(-1)^a\chi(a)
 =-\chi(2)\chi_4(d)S_d(\chi).
\]
This proves \eqref{eq:Lidentity}.

If $\chi$ is primitive and nonprincipal, then $\eta$ is odd and
primitive of conductor $4d$, since the conductors are coprime, so
$L(0,\eta)\neq0$ by \eqref{eq:Lnonzero}.
\end{proof}

\begin{lemma}\label{lem:primepowerchars}
If $d=\ell^n$ with $\ell\equiv1\pmod4$ prime, then $S_d(1)=0$,
where $1$ denotes the principal character modulo $d$. For every
nonprincipal even character $\chi$ modulo $d$, one has $S_d(\chi)\neq0$.
\end{lemma}

\begin{proof}
There are $\phi(d)/2$ units in $[1,(d-1)/2]$. The number of even
ones is
\[
 \frac{d-1}{4}-\left\lfloor\frac{d-1}{4\ell}\right\rfloor
 =\frac{\ell^n-1}{4}-\frac{\ell^{n-1}-1}{4}
 =\frac{\phi(d)}4,
\]
so $S_d(1)=0$.

If $\chi$ is nonprincipal, its inducing primitive character $\chi^*$
has conductor $\ell^a$ for some $1\leq a\leq n$. These two moduli
have the same prime support, so $\chi$ and $\chi^*$ agree as arithmetic
functions, including on nonunits. Consequently
\[
 L(s,\chi\chi_4)=L(s,\chi^*\chi_4).
\]
The latter has nonzero value at zero by Lemma~\ref{lem:Lidentity}.
Applying \eqref{eq:Lidentity} to $\chi$ proves the claim.
\end{proof}

\begin{proposition}\label{prop:primepowercriterion}
Let $\ell\equiv1\pmod4$ be prime, $n\geq1$, and let $p\neq\ell$
be a prime with $p\equiv1\pmod4$. The reduction of $B_n$ is
supersingular if and only if $p$ generates $(\Z/\ell^n\Z)^\times$.
\end{proposition}

\begin{proof}
Put $d=\ell^n$ and $H=\langle\pm p\rangle\leq U_d$. If $H$ is
proper, a nonprincipal character of $U_d/H$ gives a nonprincipal even
character modulo $d$, trivial on $H$. Supersingularity would force its
sum $S_d$ to vanish by Proposition~\ref{prop:balance}, contrary to
Lemma~\ref{lem:primepowerchars}. Conversely, if $H=U_d$, the only
character trivial on $H$ is the principal one, whose sum vanishes by
the same lemma. Thus supersingularity is equivalent to $H=U_d$.

The group $G=(\Z/d\Z)^\times$ is cyclic of order
$N=\ell^{n-1}(\ell-1)$, divisible by $4$. Write $p=g^k$ for a
generator $g$ of $G$. Its image generates $G/\{\pm1\}$ if and only
if $\gcd(k,N/2)=1$. Since $N$ and $N/2$ have the same prime divisors,
this is equivalent to $\gcd(k,N)=1$, as required.
\end{proof}

The sufficiency in this proposition is the primitive-factor result of
\cite[Proposition~4.2.6]{TT26}, in the case $\ell\equiv1\pmod4$;
the Fourier argument supplies the converse.

\begin{proposition}\label{prop:shape}
Let $d\geq3$ be odd and let $p\nmid2d$ be a prime with $p\equiv1\pmod4$.
If $J_d$ is supersingular, then $d=\ell^n$ for a prime
$\ell\equiv1\pmod4$. Moreover, $p$ generates $(\Z/\ell^n\Z)^\times$.
\end{proposition}

\begin{proof}
For each prime $\ell\mid d$, the morphism \eqref{eq:covering} makes
$J_\ell$ an isogeny factor of $J_d$, so it is supersingular.
If $\ell\equiv3\pmod4$, the group $U_\ell$ has odd order. Its
subgroup $H$ then has odd order, and no slope in \eqref{eq:CMslope}
can equal $1/2$. This is the supersingularity obstruction already
contained in the proof of \cite[Proposition~3.1.1]{TT26}.
Hence each $\ell\mid d$ is $1$ modulo $4$, and
Proposition~\ref{prop:primepowercriterion}, with $n=1$, shows that
$p$ is primitive modulo $\ell$.

Suppose that two distinct primes $\ell_1,\ell_2$ divide $d$, and put
$e=\ell_1\ell_2$. Then $J_e$ is also supersingular, and so is its
factor $A_e$. Consider the product of quadratic characters
\[
 \chi(a)=\left(\frac a{\ell_1}\right)
         \left(\frac a{\ell_2}\right).
\]
It is primitive modulo $e$, nonprincipal and even. Since a primitive
root modulo either $\ell_i$ is a nonsquare,
\[
 \chi(p)=(-1)(-1)=1.
\]
Thus $\chi$ is trivial on $\langle\pm p\rangle\leq U_e$.
Supersingularity of $A_e$ and Proposition~\ref{prop:balance} give
$S_e(\chi)=0$, contradicting Lemma~\ref{lem:Lidentity}.

It follows that $d=\ell^n$ with $\ell\equiv1\pmod4$. The primitive
factor $B_n$ is supersingular, so
Proposition~\ref{prop:primepowercriterion} gives the final assertion.
\end{proof}

\section{The Frobenius polynomial}\label{sec:frob}

In this section, let $\ell\equiv1\pmod4$ and $p\equiv1\pmod4$ be
distinct primes, and suppose that $p$ generates
$(\Z/\ell^n\Z)^\times$. Put
\[
 d_j=\ell^j,\qquad m_j=\phi(d_j)=\ell^{j-1}(\ell-1),
 \qquad h_j=m_j/2.
\]
The integer $h_j$ is even. Since reduction
$(\Z/\ell^n\Z)^\times\to(\Z/\ell^j\Z)^\times$ is surjective, the
image of $p$ generates $(\Z/\ell^j\Z)^\times$ for every $j\leq n$. Hence
Proposition~\ref{prop:primepowercriterion} shows that each $B_j$ is
supersingular, and the image of $p$ in $U_{d_j}$ has order $h_j$.

\begin{lemma}\label{lem:sparse}
For $1\leq j\leq n$, there is an integer $a_j$ such that
\begin{equation}\label{eq:sparse}
 P_{B_j,p}(X)=X^{2h_j}+a_jX^{h_j}+p^{h_j}.
\end{equation}
\end{lemma}

\begin{proof}
By Proposition~\ref{prop:CM}, the field $K_{d_j}$ acts faithfully on
$B_j$ in characteristic zero. Betti--$\lambda$-adic comparison, followed by
smooth proper base change at $p$, identifies the rational Tate module of the
reduction, compatibly with this action, with
\[
 H_1(B_j(\C),\Q)\otimes_\Q\Q_\lambda.
\]
The $\Q$-vector spaces $K_{d_j}$ and $H_1(B_j(\C),\Q)$ both have dimension
$2h_j$. Since $K_{d_j}$ is a field and acts faithfully,
$H_1(B_j(\C),\Q)$ is one dimensional over $K_{d_j}$. Hence
$V_\lambda(B_j)$ is free of rank one over
$K_{d_j}\otimes\Q_\lambda$. After extending scalars to
$\overline\Q_\lambda$, each embedding eigenspace $V_\tau$ is consequently
one dimensional.

Relation \eqref{eq:CMtwist} gives
\[
 F_p(V_\tau)=V_{\tau\circ\sigma_p^{-1}},
\]
since $uF_pv=F_p\sigma_p^{-1}(u)v$ for $v\in V_\tau$.
The element $\sigma_p=(1,\pm p)$ has order $h_j$ and fixes the
first coordinate. Hence the embeddings form two Frobenius orbits,
each of length $h_j$.

On one orbit, Frobenius sends a basis vector to a nonzero multiple of
the next, cyclically. Its characteristic polynomial is $X^{h_j}-c$,
where $c$ is the product of these multiples. The other orbit gives
$X^{h_j}-c'$. Thus
\[
 P_{B_j,p}(X)=(X^{h_j}-c)(X^{h_j}-c').
\]
Integrality of the characteristic polynomial and its constant term
$p^{h_j}$ prove \eqref{eq:sparse}.
\end{proof}

\begin{lemma}\label{lem:torsion}
For $1\leq j\leq n$,
\begin{equation}\label{eq:torsionmodule}
 B_j[2](\overline\F_p)\simeq\F_2[t]/(t^{m_j}-1),
\end{equation}
where $t$ acts as $F_p$.
\end{lemma}

\begin{proof}
The roots of $\varphi_{\ell^j}$ are obtained from
$\zeta^{2\ell^j}=-1$ by $x=\zeta+\zeta^{-1}$. They consist of $0$
and the disjoint sets
\[
 R_k=\{\zeta+\zeta^{-1}:\ord(\zeta)=4\ell^k\},
 \qquad 1\leq k\leq j.
\]
Only $\zeta$ and $\zeta^{-1}$ give the same root, so $|R_k|=m_k$.
On geometric points, $F_p$ acts by raising coordinates to the $p$-th
power. Thus a root in $R_k$ is fixed by $F_p^s$ precisely when
$p^s\equiv\pm1\pmod{4\ell^k}$. The negative sign is impossible
modulo $4$. Primitivity of $p$ modulo $\ell^k$ shows that the least
such $s$ is $m_k$. Thus $R_k$ is one Frobenius cycle.

Apply the geometric $2$-torsion description
\eqref{eq:torsionrelation} over $\overline{\F}_p$, and eliminate the generator
at $0$. This identifies the $2$-torsion with the permutation space on
the nonzero roots: the relation has coefficient $1$ at $0$, and both
the relation and the point $0$ are Frobenius-invariant. Therefore
\[
 J_{\ell^j}[2]\simeq
 \bigoplus_{k=1}^j\F_2[t]/(t^{m_k}-1).
\]
For $j=1$ this is the assertion. For $j\geq2$, use
\eqref{eq:torsionsplit} and cancel the summands of $J_{\ell^{j-1}}[2]$.
Cancellation is valid by the uniqueness of the elementary-divisor
decomposition for finite torsion modules over the PID $\F_2[t]$.
This gives \eqref{eq:torsionmodule}.
\end{proof}

\begin{lemma}\label{lem:middle}
Let $p$ be an odd prime and $A/\F_p$ a supersingular abelian variety
of positive even dimension $h$. Suppose that
\[
 P_{A,p}(X)=X^{2h}+aX^h+p^h,
 \qquad
 A[2](\overline\F_p)\simeq\F_2[t]/(t^{2h}-1)
\]
as a Frobenius module. Then $a=0$.
\end{lemma}

\begin{proof}
Fix an extension to $\overline{\Q}$ of the normalized $p$-adic valuation
$v_p$ with $v_p(p)=1$. All Frobenius roots then have valuation $1/2$.
Every monomial in the elementary symmetric expression for $a$ has
valuation $h/2$. Hence $v_p(a)\geq h/2$, and, since $h$ is even,
$a=c p^{h/2}$ for some $c\in\Z$. The roots of
$Y^2+aY+p^h$ are $h$-th powers of Frobenius roots and have complex
absolute value $p^{h/2}$. Dividing by that value shows that the roots
of $U^2+cU+1$ have absolute value $1$. Thus $|c|\leq2$.

The characteristic polynomial on the given $2$-torsion module is
$X^{2h}-1$. Reduction of $P_{A,p}$ modulo $2$ therefore makes $a$,
and hence $c$, even. We are left with $c\in\{-2,0,2\}$.

If $c=\pm2$, then
\[
 P_{A,p}(X)=(X^h\pm p^{h/2})^2.
\]
Frobenius acts semisimply on $V_2(A)$, and the binomial
$X^h\pm p^{h/2}$ is separable in characteristic zero. Thus that
binomial already annihilates $F_p$ on $V_2(A)$. By faithfulness,
\[
 F_p^h=[\mp p^{h/2}]
\]
as an endomorphism of $A$. Since the scalar is odd, $F_p^h$ would act
as the identity on $A[2]$. This contradicts the specified module:
multiplication by $t^h-1$ is nonzero modulo $t^{2h}-1=(t^h-1)^2$.
Therefore $c=0$.
\end{proof}

\begin{proof}[Proof of Theorem~\ref{thm:frobenius}]
The factors $B_j$ are supersingular, and each $h_j$ is even.
Lemmas~\ref{lem:sparse}--\ref{lem:middle} therefore give
\[
 P_{B_j,p}(X)=X^{m_j}+p^{m_j/2}.
\]
Iterating Lemma~\ref{lem:split} yields the $\F_p$-isogeny
$J_{\ell^n}\sim\prod_{j=1}^nB_j$, so the characteristic polynomials
multiply to give \eqref{eq:mainpoly}.

For a root $\alpha$ of the $j$-th factor,
$\alpha^{m_j}=-p^{m_j/2}$. Since $m_n/m_j=\ell^{n-j}$ is odd,
\begin{equation}\label{eq:commonpower}
 \alpha^{m_n}=-p^{m_n/2}
\end{equation}
for every Frobenius root of the full Jacobian. Thus, for $r=m_nb$,
\eqref{eq:extremal} gives maximality when $b$ is odd and minimality
when $b$ is even.

Conversely, the roots of the top factor are $\alpha\xi$ with
$\xi^{m_n}=1$. If the curve is maximal or minimal over $\F_{p^r}$,
all these roots have the same $r$-th power. Hence $\xi^r=1$ for every
$m_n$-th root of unity, so $m_n\mid r$. Formula
\eqref{eq:commonpower} then determines the sign and gives exactly
the stated extension degrees.
\end{proof}

\begin{proof}[Proof of Theorem~\ref{thm:classification}]
Maximality implies supersingularity by \eqref{eq:extremal}, so
(2)$\Rightarrow$(1). Proposition~\ref{prop:shape} proves
(1)$\Rightarrow$(3). Under (3), Theorem~\ref{thm:frobenius} gives
maximality over $\F_{p^{m_n}}$, proving (3)$\Rightarrow$(2).
\end{proof}

\begin{corollary}\label{cor:counts}
Under the hypotheses of Theorem~\ref{thm:frobenius}, for every $r\geq1$,
\[
 \#\Hh_{\ell^n}(\F_{p^r})
 =p^r+1-
 \sum_{\substack{1\leq j\leq n\\m_j\mid r}}
 m_j(-1)^{r/m_j}p^{r/2}.
\]
and $J_{\ell^n}(\F_p)[2^\infty]\simeq(\Z/2\Z)^n$.
\end{corollary}

\begin{proof}
The roots of the $j$-th factor are $\alpha\xi$ with $\xi^{m_j}=1$, and
their $r$-th powers sum to $m_j(-1)^{r/m_j}p^{r/2}$ if $m_j\mid r$ and to
zero otherwise; \eqref{eq:pointcount} gives the first assertion. Next,
$\#J_{\ell^n}(\F_p)=\prod_j(1+p^{m_j/2})$, and each factor is $2$ modulo
$4$ because $m_j/2$ is even. By Lemma~\ref{lem:torsion}, each $B_j[2]$ has
a one-dimensional $F_p$-fixed subspace, so $J_{\ell^n}(\F_p)[2]$ already
has order $2^n$.
\end{proof}

\begin{remark}
For $\ell=5$ and $p=13$, Theorem~\ref{thm:frobenius} gives
$L_{\Hh_{25},13}(T)=(1+13^2T^4)(1+13^{10}T^{20})$, where
$L_{C,p}(T)=T^{2g}P_{C,p}(T^{-1})$, and similarly for $\Hh_{125}$; these
recover the computations in \cite[Example~4.2.8]{TT26}.

The necessity direction is already visible at the next type of examples. For
$d=25$ and $p=157$, the prime $157$ is a primitive root modulo $5$ but not
modulo $25$ (its order modulo $25$ is $4$), so
Proposition~\ref{prop:primepowercriterion} excludes supersingularity of the
primitive factor $B_2$ and hence excludes maximality of $\Hh_{25}$ over every
finite extension of $\F_{157}$. Likewise, if $d$ has two distinct prime
factors, Proposition~\ref{prop:shape} excludes supersingularity of $J_d$ in
characteristic $p\equiv1\pmod4$; for example this applies to $d=65$ and
$p=37$.
\end{remark}

\begin{remark}\label{rem:cover}
Let $d\geq3$ be odd and $p\nmid2d$ with $p\equiv1\pmod4$. On the model
$w^2=t(t^{2d}+1)$ of $X_d$, the map
\[
 (t,w)\longmapsto(u,v)=\bigl(t^d,\;t^{(d-1)/2}w\bigr)
\]
is a morphism defined over $\F_p$ onto the elliptic curve
$E:v^2=u^3+u$, since $v^2=t^{d-1}w^2=u(u^2+1)$. As $E$ is ordinary for
$p\equiv1\pmod4$ \cite[Chapter~V, \S4]{Silverman}, $\Jac(X_d)$ is not
supersingular, although under the hypotheses of
Theorem~\ref{thm:frobenius} its quotient $\Jac(\Hh_{\ell^n})$ is.
\end{remark}

\section{The case of even degree}\label{sec:even}

For $m\geq1$, let $\mathcal C_m$ denote the smooth projective curve with
affine equation
\[
 \mathcal C_m:\qquad y^2=(x+2)\varphi_m(x).
\]
We first recall the following consequence of
\cite[Proposition~3.9]{DTT25}.

\begin{proposition}\label{prop:even-splitting}
For every $m\geq1$ and every prime $p\nmid2m$, there is an
$\F_p$-isogeny
\[
 \Jac(\Hh_{2m})\sim_{\F_p}\Jac(\Hh_m)\times\Jac(\mathcal C_m).
\]
Consequently, for every power $q$ of $p$,
\[
 \Hh_{2m}/\F_{q^2}\text{ is maximal}
 \quad\Longleftrightarrow\quad
 \Hh_m/\F_{q^2}\text{ and }\mathcal C_m/\F_{q^2}
 \text{ are both maximal}.
\]
\end{proposition}

For odd $m$, Theorem~1.4 of \cite{TT19} gives
\begin{equation}\label{eq:aux-even}
 \mathcal C_m/\F_{q^2}\text{ is maximal}
 \quad\Longleftrightarrow\quad q\equiv-1\pmod{2m}.
\end{equation}

\begin{proof}[Proof of Theorem~\ref{thm:even-p1}]
Suppose that $\Hh_d$ is maximal over $\F_{p^N}$. Since $d\geq4$, the
curve has positive genus, so $N$ is even. Write $N=2r$ and $q=p^r$.
Then $q\equiv1\pmod4$.

First assume that $4\mid d$. If $d>4$, Proposition~3.7 of \cite{DTT25}
shows that maximality implies $q\equiv-1\pmod8$, which is impossible.
Thus $d=4$. By Remark~3.8 of \cite{DTT25},
\[
 \Hh_4/\F_{q^2}\text{ is maximal}
 \quad\Longleftrightarrow\quad q\equiv-1\text{ or }5\pmod8.
\]
Since $q\equiv1\pmod4$, this is equivalent to $q\equiv5\pmod8$,
which in turn is equivalent to $p\equiv5\pmod8$ and $r$ odd. This
proves the first case.

It remains to consider $d\equiv2\pmod4$. Write $d=2m$ with $m\geq3$
odd. Proposition~\ref{prop:even-splitting} and \eqref{eq:aux-even} give
\begin{equation}\label{eq:even-reduction}
 \Hh_{2m}/\F_{q^2}\text{ is maximal}
 \quad\Longleftrightarrow\quad
 \begin{cases}
  \Hh_m/\F_{q^2}\text{ is maximal},\\
  q\equiv-1\pmod{2m}.
 \end{cases}
\end{equation}
If $\Hh_{2m}$ is maximal, then so is $\Hh_m$. By
Theorem~\ref{thm:classification},
\[
 m=\ell^a,\qquad \ell\equiv1\pmod4,
 \qquad \ord_{\ell^a}(p)=M=\ell^{a-1}(\ell-1).
\]
Theorem~\ref{thm:frobenius} also shows that $N=Mb$ with $b$ positive
and odd. This proves the necessity and shows at the same time that no
other maximal extensions can occur.

Conversely, assume these conditions, let $b$ be positive and odd, and
put $N=Mb$, $r=N/2$ and $q=p^r$. By Theorem~\ref{thm:frobenius},
$\Hh_m$ is maximal over $\F_{p^{Mb}}=\F_{q^2}$. Since $p$ generates
the cyclic group $(\Z/m\Z)^\times$ of order $M$, we have
\[
 p^{M/2}\equiv-1\pmod m.
\]
Hence
\[
 q=p^{(M/2)b}\equiv-1\pmod m.
\]
Since both $q$ and $m$ are odd, in fact $q\equiv-1\pmod{2m}$.
Equation \eqref{eq:even-reduction} now proves that $\Hh_{2m}$ is
maximal over $\F_{p^{Mb}}$.
\end{proof}

We finish by collecting the results of the present paper and the earlier
papers in one table. In all rows, $b$ denotes a positive odd integer. In
the first two rows $M=\ell^{a-1}(\ell-1)$, and in the last two rows $2k$
is the indicated multiplicative order. In these two rows $k$ is odd: the
element $p^k$ of order $2$ in $\langle p\rangle$ is $-1$ or $1+2d$, hence
$3$ modulo $4$ (the modulus is divisible by $4$), while $p\equiv3\pmod4$.

\begin{table}[ht]
\centering
\small
\renewcommand{\arraystretch}{1.22}
\begin{tabularx}{\textwidth}{@{}c c >{\raggedright\arraybackslash}X c@{}}
\toprule
$p\bmod4$ & degree & condition & maximal fields \\
\midrule
$1$ & $d$ odd
& $d=\ell^a$, $\ell\equiv1\pmod4$ prime, and $p$ generates
  $(\Z/\ell^a\Z)^\times$
& $\F_{p^{Mb}}$ \\
$1$ & $d\equiv2\pmod4$
& $d=2\ell^a$, $\ell\equiv1\pmod4$ prime, and $p$ generates
  $(\Z/\ell^a\Z)^\times$
& $\F_{p^{Mb}}$ \\
$1$ & $4\mid d$
& $d=4$ and $p\equiv5\pmod8$
& $\F_{p^{2b}}$ \\
$3$ & $d$ odd
& $\langle p\rangle\leq(\Z/4d\Z)^\times$ contains $-1$ or $1+2d$;
  $2k=\ord_{4d}(p)$
& $\F_{p^{2kb}}$ \\
$3$ & $d$ even
& $\langle p\rangle\leq(\Z/2d\Z)^\times$ contains $-1$;
  $2k=\ord_{2d}(p)$
& $\F_{p^{2kb}}$ \\
\bottomrule
\end{tabularx}
\caption{The pairs $(d,p)$ for which $\Hh_d$ is maximal over a finite
extension of $\F_p$, for $d\geq3$ and $p\nmid2d$.}
\label{tab:complete}
\end{table}

\begin{corollary}\label{cor:all-d}
Let $d\geq3$ and let $p\nmid2d$ be prime. Then $\Hh_d$ is maximal over
some finite extension of $\F_p$ if and only if one of the five cases in
Table~\ref{tab:complete} occurs. The last column of the table gives all
finite extensions over which it is maximal.
\end{corollary}

\begin{proof}
For $p\equiv1\pmod4$, the first three rows are
Theorems~\ref{thm:classification}, \ref{thm:frobenius} and
\ref{thm:even-p1}. For $p\equiv3\pmod4$ and $d$ odd, the statement is
\cite[Proposition~4.1.1]{TT26}. The case $p\equiv3\pmod4$ and $d$ even
follows from \cite[Theorem~1.5]{TT19} together with
\cite[Proposition~3.7 and Remark~3.8]{DTT25}.
\end{proof}

\ifanon\else
\section*{Acknowledgments}
S.~Tafazolian was partially supported by CNPq grant No.~302774/2025-4
and FAPESP grant No.~2024/00923-6.
\fi

\end{document}